\documentclass{amsart}
\usepackage{amssymb,latexsym}
\begin{document}
\title{Some Observations Regarding the Holomorphs of Finite Abelian Groups}

\author{Walter Becker}
\address{266 Brian Drive, Warwick, Rhode Island 02886}
\maketitle
\begin{abstract}

		Presentations for the holomorphs of abelian groups of the form $C_{p^n} \times 1^{m}$ for $p$=2 or an odd prime are given. These presentations extend the results given in Burnside's well-known text on finite groups on the holomorphs for the cyclic groups of orders $p^n$ for $p$ being an odd or even prime.
\end{abstract}

	The following set of observations all deal with the holomorphs of finite abelian groups.
The motivation here is two fold. In an extensive study of the automorphusm groups of finite groups of "low order" one frequently finds the holomorphs of certain $p$-groups arising as factors in these automorphism groups \cite{1}. A second motivation was to try an find an extension of an old result found in Burnsides' book on finite groups. In the classic text of W. Burnside \cite{2} one finds a presentation for the groups $Hol(C_t),  t= p^n$ for  all primes $p$ and for any integer $n$.  The following observations deal with extensions of this problem to other related abelian $p$-groups. When the $p$-groups are composed of different primes the problem can be decomposed into a direct product of those cases involving just one prime $p$ ,viz:
\begin{equation}
		Hol(A_p \times  B_q)  = Hol(A_p) \times  Hol(B_q).
\end{equation}
Therefore the question comes down to getting presentations for the holomorphs of  abelian $p$-groups. The following set of observations deals with a special case of abelian $p$-groups, namely when the group is a direct product of  the form
\begin{equation}
		          C_y  \times  C_p \times  C_p \times  C_p .........  where \quad y = p^n  .
\end{equation}
The cases are split up into various cases with different presentations. \\

		To generate the presentations below the use of a computer programming  system like 
CAYLEY, MAGMA or GAP will probably be required. In view of the fact that GAP appears to be exceptionally slow using a presentation input for the groups the former two programming systems would appear to be best for making use of the presentations given here \cite{3}. If one converts these presentations to a permutation presentation for the groups in question then GAP should also work very well. The other aspects of these observations, e.g., formal proofs of the statements given below, can probably be done by more standard pencil and paper theorem proof techniques. In fact the point of this set of comments is to illustrate the advantages and shortcomings of both computer and standard methods of approach to the study of finite groups and especially automorphism groups. If one computes the orders of many of the higher-order holomorphs below one readily finds that brute force computation of these groups becomes unfeasible rather quickly.\\

$\S 1.$  For the holomorph $Hol(C_y)$,  where $y = p^n$  with $n > 1$, there exists a normal subgroup of $Hol(C_y)$ isomorphic to $C_p \times  C_p$  such that $Hol(C_y)/(C_p \times  C_p)$ is isomorphic to $Hol(C_t)$  where $t=p^{(n-1)} $.\\

$\S  2.$ The automorphism group of the group $C_y \times  C_2$   ($ y = 2^n$  with $n > 1$)  is isomorphic to $( D_4 Y C_{n2} ) \times  C_2$  whose presentation can be written as:
\begin{gather}
	c^4=d^2=c^d*c=e^{n2}=c^2*e^{n4}=(c,e)=(d,e)\\ =f^2=(c,f)=(d,f)=(e,f)=1 \notag
\end{gather}
where $n2=2^{(n-2)}$ , and $n4 = \frac{n2}{2}$ , and $c^d =  d^{-1}*c*d $.\\

$\S  3.$ The action of the generators of the automorphism group of the group $C_y \times  C_2$  ( $y = 2^n$ with $n > 1$)  given in  $\S 2$ on the group $C_y \times  C_2$ is given in Table 1.\\

\begin{tabular}{|c|c|c|c|c|} \hline
\multicolumn{5}{|c|}{Table 1} \\ \hline
    &   c & d & e & f \\ \hline
 a &$ a^{-1}*b^{-1} $& $a^{-1}*b $&$a^{-5}$ & $ a^{-1}$ \\ \hline
 b & $a^{n3}*b^{-1}$ & b & b & b \\ \hline   
\end{tabular}
\linebreak

The entries in the above table are to be read in the form $a^c=a^{-1}*b^{-1}$, and in a like manner for the other entries in the table. 
Hence we have the following presentation for the $Hol(C_y \times  C_2)$  in the form
$(C_y \times  C_2) @ Aut(C_y \times  C_2)$

\begin{gather}
              a^{n1}=b^2=(a,b)=c^4=d^2=c^d*c=e^{n2}=c^2 *e^{n4}=(c,e)=(d,e)= \\ \notag
f^2=(c,f)=(d,f)=(e,f)=
	a^c*a*b=b^c*b*a^{-n3}\\ =a^d*a*b=(b,d)=a^e*a^5=(b,e)=a^f*a=(b,f)=1 \notag
\end{gather}
where $y = 2^n=n1$, and $n3=\frac{n1}{2}$. Here $n2$ and $n4$ are the same as in $\S  2.$\\

$\S 4.$ An alternative structure for the holomorph of $(C_x \times  C_2)$   ($x = 2^n$  with $n > 1$)  is given  by:
\begin{equation}
		\left[ (C_2 \times  C_2) \times  (C_2 \times  C_2) \right] @ \left[ Hol(C_2) \times  Hol(C_t) \right] .
\end{equation}

A presentation for the holomorph when written in this form is given in Table 2.\\

\begin{tabular}{|c|c|} \hline
\multicolumn{2}{|c|}{ Table 2} \\ \hline
 presentation & group \\ \hline
$a^2=b^2=c^2=d^2=$ &     $1^2\times 1^2$ \\
$=(a,b)=(a,c)=(a,d)=(b,c)=(b,d)=(c,d)=$&   \\ \hline
$e^x=f^y=g^2=(f,g)=e^f * e^{(-5)} =e^g*e$ &  $Hol(C_t)$ \\ \hline
$a*e^v=c*f^w=  $ &     actions of \\
$(a,g)=(b,g)=(c,g)=(d,g)=$ &  $hol(C_t)$ \\
$(a,e)=(b,e)=c^e*a*c=d^e*bd=$ &  on $1^2\times 1^2$ \\
$(a,f)=(b,f)=(c,f)=(d,f)= $ &          \\ \hline
 $  h^2 =$  & $Hol(C_2)$  \\ \hline
  &  actions of \\
$(a,h)=b^h*a*b=(c,h)=d^h *c*d=$ & $Hol(C_2)$ on \\ 
                &     $1^2\times 1^2$ \\ \hline
            & action of $Hol(C_2)$ \\
 $(e,h)=(f,h)=(g,h)=1$ & on $Hol(C_t)$ \\ \hline
\multicolumn{2}{|c|}{where  $t=2^{(n-1)} , x=2^n  , y=2^{(n-2)} , v=\frac{x}{2}$ and $w=\frac{y}{2}$ . }\\ \hline
\end{tabular}
\linebreak\\

Note in $\S 3$ {\itshape @ stands for a semi-direct product, but in this case the product is not a semi-direct product. The same remark applies to $\S  5$, $\S 8$, $\S 9$ and $\S 10$ below.} If we set $x = 2^{(n-1)}$ in this presentation we get the presentation for the automorphism group of the group ($ QD_{(2^n)} \times  C_2$ ). To obtain this presentation a computer programming system such as CAYLEY, GAP or MAGMA may be needed .\\

$\S  5.$ The form given in $\S 4$ for the presentations of the holomorphs of $C_x \times  C_2$ can be generalized to the cases $C_x \times  C_2 \times  C_2 \times  C_2 ..... = C_x \times 1^n$   ($ x = 2^m $) . The form in question here is:
\begin{equation}
	\left[ 1^{(n+1)} \times  (1^{(n+1)}\right]  @  \left[ Hol(1^n) \times  Hol(C_t) \right] \quad 
(  t = \frac{x}{2} )
\end{equation}
where $1^n$ means the elementary abelian group of order $2^n$, $1^{(n+1)}$ means the elementary abelian group of order $2^{(n+1)}$ and $Hol(1^n)$ means the holomorph of the elementary abelian group of order $2^n$. [ The special cases when $n$= 2 or 3 can be explicitly written out and take the form displayed in Table 3.\\

\begin{tabular}{|c|c|c|} \hline
\multicolumn{3}{|c|}{ Table 3} \\ \hline
	n=2: &  & groups\\ \hline
    &  $a^2=b^2=c^2=d^2=e^2=f^2=$&$( 1^{(n+1)} \times 1^{(n+1)} )$ \\
    &  $(a,b)=(a,c)=(a,d)=(a,e)=(a,f)= $&       \\ 
    &  $  (b,c)=(b,d)=(b,e)=(b,f)=(c,d)=$  & \\
    &   $(c,e)=(c,f)=(d,e)=(d,f)=(e,f)= $ & \\  \hline
   & $x^4=y^4=x^2*y*x^{(-2)}*y=$ &( $Hol(C_2 \times  C_2 ) = S_4 $) \\
   &  $x*y*x*(y*x*y)^{(-1)} = $ &   \\ \hline
      &$    (a,x)=b^x*a*b=c^x*a*b*c= $ &  Actions of $S_4$\\
      & $(a,y)=b^y*a*b*c=c^y*a*c=$ &   on \\
     &   $     (d,x)=e^x*d*e=f^x*d*e*f= $ & $1^3 \times 1^3 $\\
     & $(d,y)=e^y*d*e*f=f^y*d*f= $ &     group \\ \hline       
     & $s^{n1}=t^{n2}=u^2=(t,u)=s^t*s^{(-5)} = $ &  \\
     & $s^u*s= s^{n3}*a=(t^{n4}*d)= $&  $ ( Hol(C_t)$  part ) \\ \hline 
     & $(a,s)=(b,s)=(c,s)=d^s*a*d=$&  Actions of \\
     & $e^s*b*e=f^s*c*f=(s,x)=(s,y)=$ &     $Hol(C_t)$\\
     &$(a,t)=(b,t)=(c,t)=(d,t)=$ &   on the\\
     & $ (e,t)=(f,t)=(x,t)=(y,t)= $&   group\\
     & $ (a,u)=(b,u)=(c,u)=(d,u)=$& $1^3 \times 1^3$ \\
     & $(e,u)=(f,u)=(x,u)=(y,u)=1$ & \\ \hline
\end{tabular}
\linebreak

where $n1=2^m ,  n2=2^{(m-2)} , n3=\frac{n1}{2}$  and $n4=\frac{n2}{2}$. The relation in ( ) involving t, i.e., $t^{n4}*d$, is present for the cases $m  > 3$. Omitted from this table are the relations between the group $S_4$ and the group $Hol(C_t)$. Each generator of the group $S_4$ commutes with each generator of the group $Hol(C_t)$.\\

           For the case of $n$=3 the following description will enable the reader to write out the presentation in this case. Each elementary abelian group of order $2^4$ is acted upon by the group $Hol(1^3)$ in the same way and this action can be read off from the following matrix representation of the group $Hol(1^3)$ given in terms of the following four $4 \times  4 $  $(0,1)$ matrices:
\begin{gather}
t_1 = \left(\begin{matrix}
         1& 0 & 0 & 0 \\
         0 & 1 & 0 & 0 \\
         0 & 0 & 1 & 0 \\
         1 & 0 & 0 & 1
        \end{matrix}
      \right)
t_2 = \left(\begin{matrix}
         1& 0 & 0 & 0 \\
         0 & 1 & 0 & 0 \\
         0 & 0 & 1 & 0 \\
         0 & 1 & 0 & 1
        \end{matrix}
      \right)
t_3 = \left(\begin{matrix}
         1& 0 & 0 & 0 \\
         0 & 1 & 0 & 0 \\
         0 & 0 & 0 & 1 \\
         0 & 0 & 1 & 0
        \end{matrix}
      \right)    \\
t_4 = \left(\begin{matrix}
         1& 0 & 0 & 0 \\
         0 & 0 & 1 & 0 \\
         0 & 1 & 0 & 0 \\
         0 & 0 & 0 & 1
        \end{matrix}
      \right)
\end{gather}

which has the presentation:
\begin{align*}t_1^2&=t_2^2=t_3^2=t_4^2=(t_1,t_2)=(t_1,t_4)=(t_3*t_4)^3=(t_1*t_4)^4\\&=(t_1*t_3*t_2*t_3)^2=(t_2*t_3)^4=(t_2*t_4)^4=(t_2*t_3*t_4*t_3)^3\\&=t_1*t_2*t_4*t_3*t_1*t_3*t_4*t_2*t_4*t_3*t_1*t_3*t_4\\&=t_2*t_3*t_2*t_3*t_4*t_2*t_4*t_3*t_2*t_3*t_4*t_2*t_4=1.\end{align*}

This is to be read as follows: let the elementary abelian group of order $2^4$ have generators $a_1,b_1,c_1,d_1$, then $t_1$ acts on these generators as follows: $t_1$ commutes with $a_1, b_1$ and $c_1$ and for $d_1$ we have $d_1^{t_1}*a_1*d_1$. The same is true for the other three cases. The action of the group $Hol(C_t)$ on the elementary abelian group $1^8$ is analogous to the $n$=2 case with the s generator commutating with the first four generators and then on the next four generators $s^{n3}*a_1=a_2^s*a_1*a_2=b_2^s*b_1*b_2=c_2^s * c_1*c_2=d_2^s*d_1*d_2=1$. The order of the holomorph of $C_2 \times  C_2 \times  C_2$ is 1344. The order of the next case, i.e., $n$=4 is 24*8!/2 = 322,560 so the orders of these factors get very large very fast. It would be nice if one could come up with an iterative scheme for determining the groups $ Hol(1^n)$ and their actions on the elementary abelian 
groups of order $2^{(n+1)}$ , in terms of  a set of $(n+1) \times  (n+1),  (0,1)$ matrices. \\

$\S 6.$ The automorphism groups of the groups $C_y \times  C_3$  where $y = 3^n$ , $n > 1$ , are 
isomorphic to  $( [ Qf(3) Y C_z ] @ C_2 ) \times  C_2$  or ($ [Qf(3) @ C_2] Y C_z ) \times  C_2$ :
\begin{gather*}          
 a^3=b^3=c^3=(a,b)=(a,c)=b^c *a^{(-1)} * b^{(-1)} =\\f^z = (a,f)=(b,f)=(c,f)=a*f ^{(-z1)}=\\
         d^2 =(a,d)=b^d * b^{(-s)} = cd * c^{(-t)} =\\ e^2 = (a,e)=(b,e)=(c,e)=d,e)=1
\end{gather*}
where $s^{(p-1)} = 1 mod(p)$  and $t = s^{(p-2)} mod(p)$. That is, $s$   is a $(p-1)$-th root of $p$. Here 
 $z = 3^{(n-1)}$  and $z1 = 3^{(n-2)} $.\\

$\S  7.$ For the presentation of the groups $Aut(C_y \times  C_3)$   given in $\S  6$  the actions of the generators of $Aut(C_y \times  C_3)$ on the generators of $C_y \times  C_3$  are given explicitly in Table 4.\\

\begin{tabular}{|c|c|c|c|c|c|c|}\hline
\multicolumn{7}{|c|}{Table 4} \\ \hline
\multicolumn{7}{|c|}{Generators} \\ \hline
 &     a     &       b  &            c    &           d      &           e     &        f  \\ \hline 
   P  &$     P^{ap} $ &          P     &     $   P^{tt} * Q$&  $    P^{-1} $&   $ P^x $&$ P^af$\\ \hline
   Q   &       Q  &        P*Q   &           Q   &              Q   &  $        Q^{-1} $&     Q \\ \hline
\end{tabular}\\
\linebreak

and therefore a presentation of the holomorph of $(C_y \times  C_3)$ is:
\begin{equation*}
        a^3=b^3=c^3=(a,b)=(a,c)=b^c *a^{(-1)} * b^{(-1)} =
\end{equation*}
\begin{equation*}
 f^z = (a,f)=(b,f)=(c,f)=a*f^{(-z1)} =
\end{equation*}
\begin{equation*}
                     d^2 =(a,d)=b^d * b = c^d * c=e^2 = 
\end{equation*}
\begin{equation*}
(a,e)=(b,e)=(c,e)=(d,e)=
\end{equation*}
\begin{equation*}
                   	 P^y=Q^3=(P,Q)=
\end{equation*}
\begin{equation*}
       P^a*P^{(-ap)} =(b,P) =P^c*P^{(-tt)} * Q^{(-1)} =P^d*P = 
\end{equation*}
\begin{equation*}
                  P^e*P^{(-x)} =P^f*P^{(-af)} =
\end{equation*}
\begin{equation*}
                    (Q,a)=Q^b*P^{(-1)} *Q^{(-1)} =(c,Q) =(d,Q)=Q^e*Q=(Q,f)=1 .
\end{equation*}

     Here  $ (tt) = 3^{(n -1)}, (ap)  = [2*3^{(n-1)} + 1] $ ,  (af) is a   $ 3^{(n-1)}$ -th root of unity and  $x = (p-1)$-th root of unity  
   (i.e.,  $x^{(p-1)}  = 1 mod(p^n )$ ).\\

$\S  8.$ The holomorph of $C_y \times  C_3$ can be written in the following alternate form similar to that for the 2-group in $\S  4$:\\

	$[C_3 \times  C_3) \times  (C_3 \times  C_3) ] @  [ Hol(C_p) \times  Hol(C_t)]$\\

where $y=3^n $ and $t=3^{(n-1)}$ .

 A presentation for this case can be written in the form broken down into various parts and is given in Table 5.\\

\begin{tabular}{|c|c|} \hline
\multicolumn{2}{|c|}{Table 5} \\ \hline
Presentation & Group \\ \hline
 & \\
 $a^3=b^3=c^3=d^3=(a,b)=(a,c)=$ & \\
$(a,d)=(b,c)=(b,c)=(c,d)=$ &$ \left(1^4\right)$ \\ 
    &   \\ \hline \hline 
  & \\                                       
   $ (e1)^3 = (f1)^2=(e1)^{(f1)} *(e1)= $ &  $ \left( Hol(C_3)\right) $ \\ 
     &  \\ \hline \hline
   & actions  \\                                                                                 
 $a^{(e1)} *a^{(-1)} *b^{(-1)} = (b,(e1))=a^{(f1)} *a = $ & of  \\
$(b,(f1))=c^{(e1)} *c^{(-1)} * d^{(-1)} = (d,(e1))=c^{(f1)} *c =$ &$ \left( Hol(C_3)\right) $   \\ 
 $(d,(f1))=$  & on $(1^4)$  \\ \hline \hline
    &    \\
$e^x =f^y=e^f*e^4 =$  & $ \left( Hol(C_t) \right)$ \\  
    &   \\   \hline \hline    
     &   \\                                                                                                    
 $e^{x2} * b=d^{(-1)}* f^{y1}= (e,(e1))=(e,(f1))=(f,(e1))=$ & actions \\
 $(f,(f1))=(a,e)= c^e *a^{(-1)} * b^{(-1)} =a^f * a = $ & of \\
$(c,f)=(b,e)=d*e*b*d^{(-1)} =b^f *b =(d,f)=1$ &$ \left( Hol(C_t) \right)$  \\ 
    &   \\ \hline \hline
\end{tabular} 
\linebreak \\

where $x=3^n  , x1=2*3^{(n-2)} , x2= 3^{(n-1)} , y=2*3^{(n-1)}$ and $y1=2*3^{(n-2)}$. Note that an alternate form occurs when one makes the replacements:\\

        \qquad   \qquad   $e^f *e^4 \qquad  \qquad $          by  \qquad \qquad $  e^f *e^{(-2)}  $ \\
and \\

 \qquad \qquad          $ d^{(-1)} * f^{y1} \quad \quad $  by \qquad \qquad$ d^{(-1)} *f^{yt} \qquad $where $yt = 4 * 3^{(n-2)}$ . \\

$ \S 9.$ The results of $\S  6$ can be generalized to any odd prime $p$ as follows:

1. The automorphism group of  $Aut(C_y \times  C_p)$ where $y = p^n$  and $n > 1$ is isomorphic to
\begin{equation}
     \left( \left[ Qf(p) Y C_z \right ] @ C_{(p-1)} \right) \times  C_{(p-1)}
\end{equation}
 or
\begin{equation}
 \left( \left[ Qf(p) @ C_{(p-1)} \right] Y C_z  \right ) \times  C_{(p-1)} 
\end{equation}
with the presentation:
\begin{gather*}          
               a^p=b^p=c^p=(a,b)=(a,c)=b^c*a^{(-1)} * b^{(-1)}\\ 
	=f^z = (a,f)=(b,f)=(c,f)=a*f^{(-z1)} =
         d^{(p-1)}\\ =(a,d)=bd * b^{(-s)} = c^d * c^{(-t)} =e^{(p-1)}\\
      = (a,e)=(b,e)=(c,e)=d,e)=1
\end{gather*}
 where $s^{(p-1)} = 1 mod(p)$  and $t = s^{(p-2)} mod(p)$. That is, s is a $(p-1)$-th root of $p$. Here 
 $z = p^{(n-1)}$ , and $z1 = p^{(n-2)}$ .\\

2. The actions of the generators of this automorphism group on the generators of $C_y \times  C_p$ is given in Table 6. \\

\begin{tabular}{|c|c|c|c|c|c|c|} \hline
\multicolumn{7}{|c|}{Table 6} \\ \hline
\multicolumn{7}{|c|}{Generators} \\ \hline
   &    a  &        b&         c &          d  &         e  &         f \\ \hline  
 P &    $ P^{(ap)}$   &      P  &      $   P^{(tt)} * Q $&  $ P^{(-1)} $& $  P^x $   &$P^{(af)}$\\ \hline
 Q &      Q  &   P*Q   &     Q  &$  Q^{(y1)}$ &    $ Q^{(y1)}$  &Q \\ \hline
\end{tabular} \\

Here   $(tt) = p^{(n -1)}$, $(ap)  = \left[(p-1)*p^{(n-1)} + 1\right]$, and $(af)$ is a \, $p^{(n-1)}$-th  root of unity . In addition we have $x^{(p-1)} = 1 mod(p^n)$  and $(y1)^{(p-1)}  = 1 mod(p)$ . That is, $x$ and $(y1)$ are also $p-1$ roots of unity in 
$GF(p^n)$  and GF(p) respectively.

3. Combining 1 and 2 will yield a presentation for the $Hol(C_y \times  C_p)$. \\

$\S 10.$  The results of $\S 8$ can be generalized, like that in $\S  5$ for the 2-groups to the cases for the groups $C_y \times  C_p \times  C_p ...... = C_y \times  1^n$ , where now the structure of the holomorphs is :

\begin{equation}
	\left[ \left(1^{(n+1)}\right) \times   \left(1^{(n+1)}\right) \right] @ \left[ Hol(1^n) \times  Hol(C_t)  \right]. 
\end{equation}
Here $y = p^m$  and $t =p^{(m-1)}$ . Note that in these cases the orders of the groups are sufficiently large that a direct calculation of these groups by computers can be rather time consuming. 
 The advantage of  writing the holomorph in this form is that if one can find a representation of the 
$Hol(1^n)$ valid for all $n$ and $p$ then one would not need to compute the automorphism groups of the groups  $( C_y \times  1^n)$ for each different case n, assuming one can explicitly find a representation for the automorphism groups of  $(C_y \times  1^n)$ that is valid for all $p$ with a fixed n.\\

$\S  11.$  From $\S  8$ and $\S 10$ one knows that the group $Hol(1^n)$  appearing in these presentations is a subgroup of $Aut( 1^{(n+1)} ) = GL(n+1 , p)$.  There is more than one conjugacy class of subgroups of $Hol( 1^n )$ in  $GL(n+1 ,p)$. For the case of $n$=1, computer results suggest that one has $(p-1)$ inequivalent classes of $Hol(C_p)$ . Can this be proven for all odd primes $p$?
For the cases when $n > 1$ the only computer runs available are for $p$=3  and $n$=2. Here one finds 
four inequivalent classes of $Hol( C_3 \times  C_3)$. Can one generalize this to the case of \\

a) $n$=2 and arbitrary odd prime $p > 3$, and \\
 
b) $n>2$ and arbitrary odd prime $p$ greater than or equal to 3 ?\\

Note that the particular choice of the $Hol( 1^n)$  that appears in the presentations for  the structure given in $\S 12$ when it acts on each of the elementary abelian groups of order $p^{(n+1)}$ yields a group with the center $C_p$. In the case of $n$=2 and $p$=3 there is only one such case. For the other cases is it also true that this case is unique ?\\

$\S 12.$ Some additional suggestion for future studies. We have given some indications of additional studies that might be interesting to pursue above. Two other curiousities that might be of interest to the reader to follow up on are:\\

a:) We have not constructed the isomorphic map between these two different representations of $Hol\left(C_x\times C_2\right)$. It might be an interesting excercise for the reader to try and construct the mapping between the presentations of $Hol\left(C_x\times C_2\right)$ given in $\S 3$ and $\S 4.$, and likewise between the forms of  $Hol(C_y \times  C_3)$ given in $\S 7$ and $\S 8$.\\

b:) One could also proceed as above with other forms, e.g., the holomorphs of the abelian groups $C_s \times C_t \times 1^n$ where $s=p^k$ and $t=p^m$ with $k>m$. The closed form expressions for these cases, if they can be found, might be rather complicated and not very useful. The interest in the cases discussed above are of interest because they arise rather naturally as the automorphism groups of certain finite groups. It is unclear if these other cases would arise in such a natural fashion, but it might be worth investigating, and the effort might yield some interesting results.\\
 
$\S 13.$ Acknowledgements.  Most of this work was done at Brown University on a DEC VMS computer in the Department of Linguistics and Cognitive Sciences. The author wishes to thank Dr. James Anderson of this Department for giving him the requisite time for performing these and a great many other calculations dealing with the Automorphism Groups of Finite Groups. The programming system used was CAYLEY which was made available to the author by the generousity of Dr. John Cannon of the Department of Mathematics,  University of Syndey, Australia.\\


\begin{thebibliography}{99}
\bibitem{3} Explaination of the slow response in GAP to a presentation input:
"The reason is very simple: in Cayley a finite fp group is represented by
a cayley graph and Cayley contained algorithms specially designed for the 
cayley graph representation of a group. Gap simply uses the regular 
permutation representation in all cases." (John Cannon private response to this question as to why GAP is so slow computating groups with a presentation input.)
\end{thebibliography}
\end{document}